\documentclass{amsart} 

\usepackage{amsmath,amsthm}     
\usepackage{amssymb}            
\usepackage{euscript}           
\usepackage{enumerate,calc}
\usepackage[matrix,arrow,curve,frame]{xy}    

\xymatrixcolsep{1.9pc}                          
\xymatrixrowsep{1.9pc}
\newdir{ >}{{}*!/-5pt/\dir{>}}                  

\newtheorem{thm}[subsection]{Theorem}
\newtheorem{defn}[subsection]{Definition}
\newtheorem{prop}[subsection]{Proposition}

\newtheorem{cor}[subsection]{Corollary}

\theoremstyle{definition}  
\newtheorem{example}[subsection]{Example}

\newtheorem{remark}[subsection]{Remark}

  {\end{list}}

\newcommand{\dfn}{\textbf} 

\newcommand{\mdfn}[1]{\dfn{\mathversion{bold}#1}} 

\newcommand{\Smash}             {\wedge}
\newcommand{\sm}             {\wedge}

\newcommand{\Wedge}             {\vee}

\newcommand{\tens}              {\otimes}               
\newcommand{\iso}               {\cong}

\newcommand{\cat}{\EuScript}    
\newcommand{\cA}{{\cat A}}      
\newcommand{\cB}{{\cat B}}      
\newcommand{\cC}{{\cat C}}

\newcommand{\cD}{{\cat D}}

\newcommand{\cG}{{\mathcal G}}

\newcommand{\cM}{{\cat M}}
\newcommand{\cN}{{\cat N}}

\newcommand{\cR}{{\cat R}}
\newcommand{\cS}{{\cat S}}
\newcommand{\cT}{{\cat T}}

\newcommand{\sSet}{s{\cat Set}}

\newcommand{\Ho}{\text{Ho}\,}
\newcommand{\ho}{\text{Ho}\,}

\newcommand{\field}[1]  {\mathbb #1} 

\newcommand{\F}         {\field F}

\newcommand{\Z}         {\field Z}
\newcommand{\bZ}         {\field Z}

\newcommand{\Q}         {\field Q}
\newcommand{\bQ}         {\field Q}

\DeclareMathOperator{\Hom}{Hom}
\DeclareMathOperator{\Fhom}{\Hom}
\DeclareMathOperator{\Fhomu}{\mM_{\Ch}}
\DeclareMathOperator{\NFhomu}{\mN_{\Ch}}

\DeclareMathOperator{\Ext}{Ext}

\DeclareMathOperator{\RHom}{RHom}

\DeclareMathOperator{\sDer}{{\mathbf{Der}}}
\DeclareMathOperator{\sHH}{{\mathbf{HH}}}
\DeclareMathOperator{\Der}{Der}
\DeclareMathOperator{\HH}{HH}
\DeclareMathOperator{\sTDer}{{\mathbf{TDer}}}
\DeclareMathOperator{\sTHH}{{\mathbf{THH}}}
\DeclareMathOperator{\TDer}{TDer}
\DeclareMathOperator{\THH}{THH}
\DeclareMathOperator{\Aut}{Aut}

\DeclareMathOperator{\Mod}{Mod-}

\DeclareMathOperator{\lmod}{-Mod}

\DeclareMathOperator{\im}{Im}

\newcommand{\ra}{\rightarrow}                   
\newcommand{\lra}{\longrightarrow}              
\newcommand{\la}{\leftarrow}                    
\newcommand{\from}{\leftarrow}                    
\newcommand{\lla}{\longleftarrow}               
\newcommand{\llra}[1]{\stackrel{#1}{\lra}}      
\newcommand{\llla}[1]{\stackrel{#1}{\lla}}      

\newcommand{\we}{\llra{\sim}}                   
\newcommand{\bwe}{\llla{\sim}}

\newcommand{\dbra}{\rightrightarrows}           

\newcommand{\blank}{-}                          
\newcommand{\und}{\underline}
\newcommand{\mM}{\underline{\cM}}

\newcommand{\adjoint}{\rightleftarrows}

\newcommand{\bdd}[1]{\partial\Delta^{#1}}
\newcommand{\del}[1]{\Delta^{#1}}

\newcommand{\he}{\simeq}

\newcommand{\Fp}{{\F_p}}
\newcommand{\Ft}{{\F_2}}

\newcommand{\PP}{{\field P}}   

\newcommand{\HZ}{H\Z}
\newcommand{\HZt}{H\Z/2}

\numberwithin{equation}{subsection}

\newenvironment{myequation}
  {\addtocounter{subsection}{1}\begin{equation}}
  {\end{equation}$\!\!$}

\newcommand{\Alg}{-\cA lg}
\newcommand{\DGA}{\cD G\cA}

\newcommand{\kDGA}{k\!-\!\cD GA}
\newcommand{\ZDGA}{\Z\!-\!\cD GA}
\newcommand{\RSp}{\cR ing\cS p}
\newcommand{\Ch}{Ch}
\newcommand{\rMod}{\Mod}

\newcommand{\Qdga}{$\Q$-dga\ }

\newcommand{\QdgasP}{$\Q$-dgas}
\newcommand{\kdgas}{k-dgas\ }

\newcommand{\ChZ}{\Ch}
\newcommand{\chk}{\Ch_k}
\newcommand{\mc}{\colon \,}

\newcommand{\Sp}{Sp}
\DeclareMathOperator{\hEnd}{hEnd}
\newcommand{\RingSpectra}{S\Alg}
\newcommand{\Algebra}{\Alg}

\newcommand{\mN}{\underline{\cN}}

\newcommand{\Ab}{\cat Ab}
\newcommand{\sab}{s{\cat Ab}}

\newcommand{\Spectra}{\Sp^\Sigma}

\newcommand{\Ring}{\cat{R}ing}

\newcommand{\namma}{\nut}
\newcommand{\nut}{\nu}
\newcommand{\xib}{\bar\xi}

\DeclareMathOperator{\ch}{ch}

\begin{document}

\title{Topological equivalences for differential graded algebras}
                                                           
\date{\today;\ \  2000 AMS Math.\ Subj.\ Class.: 55U99, 55P43, 18G55}                                                          

\author{Daniel Dugger}
\address{Department of Mathematics\\ University of Oregon\\ Eugene, OR
97403}

\email{ddugger@math.uoregon.edu}

\author{Brooke Shipley}
\thanks{Second author partially supported by NSF Grant No. 0134938 and a
Sloan Research Fellowship}

\address{Department of Mathematics\\ University of Illinois at Chicago
\\ Chicago, IL 60607}

\email{bshipley@math.uic.edu}

\begin{abstract}
We investigate the relationship between differential graded
algebras (dgas) and topological ring spectra.  Every dga $C$ gives
rise to an Eilenberg-Mac\,Lane ring spectrum denoted $HC$.
If $HC$ and $HD$ are weakly equivalent, then we say $C$ and $D$
are topologically equivalent.    Quasi-isomorphic dgas are
topologically equivalent, but we produce explicit counter-examples of
the converse.  We also develop an associated notion of topological
Morita equivalence using a homotopical version of tilting.
\end{abstract}

\maketitle

\tableofcontents

\section{Introduction}

This paper deals with the relationship between differential graded
algebras (dgas) and topological ring spectra.  DGAs are considered
only up to quasi-isomorphism, and ring spectra only up to weak
equivalence---both types of equivalence will be denoted $\he$ in what
follows.  Every dga $C$ gives rise to an Eilenberg-Mac\,Lane ring
spectrum denoted $HC$ (recalled in Section~\ref{se:EM}), and of course
if $C\he D$ then $HC \he HD$.  It is somewhat surprising that the
converse of this last statement is not true: dgas which are not
quasi-isomorphic can give rise to weakly equivalent ring spectra.  If
$HC\he HD$ we will say that the dgas $C$ and $D$ are
\dfn{topologically equivalent}.  Our goal in this paper is to
investigate this notion, with examples and applications.  The papers
~\cite[Section 3.9]{keller} and~\cite{ship-expos} give expository
accounts of some of this material.

\medskip

\subsection{Explicit examples}
\label{se:explicit}
In the first sections of the paper we are concerned with producing
examples of dgas which are topologically equivalent but not
quasi-isomorphic.  One simple example turns out to be
$C=\Z[e;de=2]/(e^4)$ and $D=\Lambda_{\F_2}(g)$ with the degrees $|e| =
1$ and $|g| = 2$.  That is to say, $C$ is a truncated polynomial
algebra with $de=2$ and $D$ is an exterior algebra with zero
differential.  To see that these dgas are indeed topologically
equivalent, we analyze Postnikov towers and their associated
$k$-invariants.  This reduces things to a question about the
comparison map between Hochschild cohomology $\HH^*$ and topological
Hochschild cohomology $\THH^*$, which one can resolve by referring to
calculations in the literature.

Similar---but more complicated---examples exist with $\F_2$ replaced
by $\F_p$.  It is an interesting feature of this subject that as $p$
grows the complexity of the examples becomes more and more intricate.

\subsection{Model categories}
\label{se:modelcat}
The above material has an interesting application to model categories.
Recall that a Quillen equivalence between two model categories
is a pair of adjoint functors satisfying certain axioms with respect
to the cofibrations, fibrations, and weak equivalences.  Two model
categories are called Quillen equivalent if they can be
connected by a zig-zag of such adjoint pairs.  Quillen equivalent
model categories represent the same `underlying homotopy theory'.

By an \dfn{additive model category} we mean one whose underlying
category is additive and where the additive structure behaves well
with respect to `higher homotopies'; a precise definition is given in
\cite[Section 2]{additive}.  So this excludes the model categories of
simplicial sets and topological spaces but includes most model
categories arising from homological algebra.  Two additive model
categories are \dfn{additively Quillen equivalent} if they can be
connected by a zig-zag of Quillen equivalences in which all the
intermediate steps are additive.  The following is a strange and
interesting fact:
\begin{itemize}
\item {\it It is possible for two additive model categories to be Quillen
equivalent but not additively Quillen equivalent.\/}
\end{itemize}
One might paraphrase this by saying that there are two
`algebraic' model categories which have the same underlying homotopy
theory, but where the equivalence cannot be seen using only algebra!
Any zig-zag of Quillen equivalences between the two must necessarily
pass through a non-additive model category.

In Section~\ref{se:QEnotAQE} we present an example demonstrating the above
possibility.  It comes directly from the two dgas $C$ and $D$ we wrote
down in (\ref{se:explicit}).  If $A$ is any dga, the category of
differential graded $A$-modules (abbreviated as just `$A$-modules'
from now on) has a model structure in which the weak equivalences are
quasi-isomorphisms and the fibrations are surjections.  We show that
the model categories of $C$-modules and $D$-modules are Quillen
equivalent but not additively Quillen equivalent.

\medskip

\subsection{Topological tilting theory for dgas}
\label{se:intro-tilting}
The above example with model categories is actually an application of
a more general theory.  One may ask the following question: Given two
dgas $A$ and $B$, when are the categories of $A$-modules and
$B$-modules Quillen equivalent?  We give a complete answer in terms of
a homotopical tilting theory for dgas, and this involves topological
equivalence.

Recall from Morita theory that two rings $R$ and $S$ have equivalent
module categories if and only if there is a finitely-generated
$S$-projective $P$ such that $P$ is a strong generator and the
endomorphism ring $\Hom_S(P,P)$ is isomorphic to $R$.  Rickard
\cite{R} developed an analogous criterion for when the derived
categories $\cD(R)$ and $\cD(S)$ are equivalent triangulated
categories, and this was extended in \cite[4.2]{dereq} 
to the model category level.  Explicitly, the model categories of chain
complexes $\Ch_R$ and $\Ch_S$ are Quillen equivalent if and only if
there is a bounded complex of finitely-generated $S$-projectives $P_*$
which is a weak generator for the derived category $\cD(S)$, and whose
endomorphism dga $\Hom_S(P,P)$ is quasi-isomorphic to $R$ (regarded as
a dga concentrated in dimension zero).

We wish to take this last result and allow $R$ and $S$ to be dgas
rather than just rings.  Almost the same theorem is true, but
topological equivalence enters the picture:

\begin{thm}
\label{th:tilting}
Let $C$ and $D$ be two dgas.  The model categories of $C$- and
$D$-modules are Quillen equivalent if and only if there is a cofibrant
and fibrant representative $P$ of a compact
generator in $\cD(C)$ whose endomorphism dga $\Hom_C(P,P)$
is topologically equivalent to $D$.
\end{thm}

See Definition~\ref{def:CompGen} for the definition of a compact generator.
If $\cD(C)$ has a compact generator, then a cofibrant
and fibrant $C$-module representing this generator always exists.

If one takes Theorem~\ref{th:tilting} and replaces `topologically
equivalent' with `quasi-isomorphic', the resulting statement is
false.  For an example, take $C$ and $D$ to be the two dgas already
mentioned in Section~\ref{se:explicit}.  The model categories of $C$-
and $D$-modules are Quillen equivalent since they only depend on $HC$
and $HD$, but it is easy to show (see
Section~\ref{se:QEnotAQE}) that no $D$-module can have $C$ as its
endomorphism dga.

We have the following parallel of the above theorem, however:

\begin{thm}\label{thm-add}
Let $C$ and $D$ be two dgas.  The model categories of $C$- and
$D$-modules are \dfn{additively} Quillen equivalent if and only if there is
a cofibrant and fibrant representative $P$ of a compact generator in $\cD(C)$
such that $\Hom_C(P,P)$ is quasi-isomorphic to $D$.
\end{thm}

The above theorem should be compared to a result of Keller involving
derived equivalences for dgas (with many objects) \cite[3.11]{keller}.  
In fact, by combining our result with Keller's one finds that 
two dgas have module categories which are
additively Quillen equivalent if and only if  the dgas are ``dg Morita
equivalent'' in the sense of \cite{keller}.  This is discussed in
detail in Section~\ref{se:dgcats}.

\subsection{Topological equivalence over fields}
We do not know a general method for deciding whether two given dgas
are topologically equivalent or not.  The examples of topological
equivalence known to us all make crucial use of dealing with dgas over
$\Z$.  It would be interesting to know if there exist nontrivial
examples of topological equivalence for dgas defined over a field.
Here is one negative result along these lines, whose proof is given in
Section~\ref{se:fields}.

\begin{prop}\label{prop-Q-same}
If $C$ and $D$ are both \QdgasP, then they are topologically
equivalent if and only if they are quasi-isomorphic.
\end{prop}

\subsection{Organization}
The main results of interest in this paper---described above---are contained in
Sections~\ref{se:examples}, \ref{se:tilting}, and \ref{se:QEnotAQE}.
Readers who are impatient can jump straight to those sections.

Sections~\ref{se:background} through \ref{se:Postnikov-ringsp}
establish background on dgas and ring spectra needed for the examples
in Section~\ref{se:examples}.  This background material includes
Postnikov sections, $k$-invariants, and the role of Hochschild and
topological Hochschild cohomologies.

The tilting theory results are given in Section~\ref{se:tilting}.  To
prove these, one needs to use certain invariants of stable model
categories---namely, the homotopy endomorphism ring spectra of
\cite{hend}.  These invariants are defined very abstractly and so are
difficult to compute, but in Section~\ref{se:hoend} we state some
auxiliary results simplifying things in the case of model categories
enriched over $\Ch_\Z$.  The proofs of these results are rather technical
and appear in~\cite{additive}.

\subsection{Notation and terminology}
If $\cM$ and $\cN$ are model categories, a Quillen pair $L\colon \cM
\adjoint \cN \colon R$ will also be referred to as a \dfn{Quillen map}
$L\colon \cM \ra \cN$.  The terms `strong monoidal-' and `weak
monoidal Quillen pair' will be used often---they are defined in
\cite[3.6]{SS3}.  
A \dfn{strong monoidal Quillen pair} is basically
a Quillen pair between monoidal model categories where $L$ is
strong monoidal.  Finally, if $\cC$ is a category we write $\cC(X,Y)$
for $\Hom_\cC(X,Y)$.
Throughout the paper, for a symmetric spectrum $X$ the notation $\pi_*X$
denotes the derived homotopy groups; that is, one first replaces $X$
by a fibrant spectrum and then evaluates $\pi_*$.

Finally, one notational convention we use throughout the paper is to
denote the degree of an element as a subscript.  Thus, $\F_2[x_3]$
denotes a polynomial ring on a class $x$ in degree $3$.  We will
sometimes drop the subscript when the notation becomes too cluttered.   


\section{Background on dgas and ring spectra}
\label{se:background}
In this section we review model category structures on dgas and ring
spectra.  We also recall the construction which associates to every dga a
corresponding Eilenberg-Mac\,Lane ring spectrum.

\medskip

\subsection{DGAs}
If $k$ is a commutative ring,
let $\Ch_k$ denote the category of (unbounded) chain complexes of
$k$-modules.  Just to be clear, we are grading things so that the
differential has the form $d\colon C_n\ra C_{n-1}$.  Recall that
$\Ch_k$ has a model category structure---called the `projective' model
structure---in which the weak equivalences are the quasi-isomorphisms
and the fibrations are the surjections \cite[2.3.11]{hovey}. 
The tensor product of chain complexes makes this into a symmetric monoidal
model category in the sense of \cite[4.2.6]{hovey}. 

By a \mdfn{$k$-dga} we mean an object $X\in \Ch_k$ together with maps
$k[0]\ra X$ and $X\tens X \ra X$ giving an associative and unital
pairing.  Here $k[0]$ is the complex consisting of a single $k$
concentrated in dimension $0$.  Note that we only require associativity
here, not commutativity.

By \cite[4.1(2)]{SS1}, 
a model category structure on \kdgas can be
lifted from the projective model structure on $\Ch_k$.
The weak equivalences of \kdgas are again the quasi-isomorphisms and
the fibrations are the surjections.  The cofibrations are then
determined by the left lifting property with respect to the acyclic
fibrations.  A more explicit description of the cofibrations comes
from recognizing this model structure as a cofibrantly generated model
category~\cite[2.1.17]{hovey}.  Let $k(S^n)$ be the free $k$-algebra
on one generator in degree $n$ with zero differential.  Let
$k(D^{n+1})$ be the free $k$-algebra on two generators $x$ and $y$
with $|x|=n$ and $|y|=n+1$ such that $dx=0$ and $dy = x$.  The
generating cofibrations are the inclusions $i_n \mc k(S^n) \to
k(D^{n+1})$ and the generating acyclic cofibrations are the maps $j_n
\mc 0 \to k(D^{n+1})$.  These generating maps are constructed from the
generating maps for $\chk$~\cite[2.3.3]{hovey} by applying the free
tensor algebra functor $T_k$.  Given $C$ in $\chk$, $T_k(C) =
\oplus_{n \geq 0} C^{\otimes n}$ where $C^{\otimes 0} = k[0]$ and all
tensor products are over $k$.

\begin{remark}
Cofibrant replacements of \kdgas play an important role in what
follows.  There is a functorial cofibrant replacement
arising from the cofibrantly generated structure~\cite{SS1}, but this
gives a very large model.  We sketch a construction of a smaller
cofibrant replacement which is useful in calculations.  Suppose given $C$ in
$\kDGA$ with $H_iC=0$ for $i<0$.
Choose generators of $H_*C$ as a $k$-algebra.
Let $G$ be the free graded $k$-module on the
given generators and let $T_1C = T_k(G)$ be the associated free
$k$-dga
with zero differential.
Define $f_1\mc T_1C \to C$
by sending each generator in $G$ to a chosen cycle representing the
associated generator in $H_*C$.  The induced map $(f_1)_*$ in homology
is surjective.  Let $n$ be the smallest degree in which $(f_1)_*$ has
a kernel, and pick a set of $k$-module generators for the kernel in
this dimension.
For each chosen generator
there is an associated cycle in $T_1C$ and thus a map
$k(S^n) \to T_1(C)$.  The pushout $k(D^{n+1}) \from k(S^n) \to T_1(C)$
has the effect of killing the associated element in $H_n[T_1(C)]$.  Let
$T_2C$ be the pushout of
\[ \coprod k(D^{n+1}) \from \coprod k(S^{n}) \to T_1C
\]
where the coproduct runs over the chosen generators.  (Note
that the coproduct is taken in the category $\kDGA$, and so is
complicated---it's analogous to an amalgamated product of groups.)
Since the elements being killed here are in the kernel of $(f_1)_*$,
each map $k(S^{n}) \to T_1(C) \to C$ extends to a map $k(D^{n+1})
\to C$.  Thus, there is a map from the pushout $f_2\colon T_2C \to C$.
The induced map $(f_2)_*$ is an isomorphism up through degree $n$,
and now one repeats the process in degree $n+1$.
The colimit of the resulting $T_iC$'s is a cofibrant replacement
$T_{\infty}C \to C$.  (Note that $T_\infty C\ra C$ is not a fibration,
however.)
\end{remark}

\begin{example}
\label{ex:cofibrep}
By way of illustration, let $C=\Z/2$ considered as a $\Z$-dga
(concentrated in degree $0$).  We may take $T_1C=\Z$.  The element $2$
must now be killed in homology, so we attach a free generator $e$ in
degree $1$ to kill it.  That is, we form the pushout of a diagram
\[ \Z(D^1) \la \Z(S^0) \ra T_1C.
\]
The result is $T_2(C)=\Z[e; de=2]$ with $|e|= 1$.  The map $T_2C \ra
\Z$ is an isomorphism on $H_0$ and $H_1$, but on $H_2$ we have a
kernel (generated by $e^2$).  After forming the appropriate pushout we
have $T_3C=\Z\langle e, f; de=2, df=e^2\rangle$ with $|e| = 1, |f| =
3$ (forgetting the differentials, this is a tensor algebra on $e$ and
$f$).

Next note that $H_3(T_3C)= 0$, but $H_4(T_3C)$ is nonzero---it's
generated by $ef+fe$.  So we adjoin an element $g$ with $dg=ef+fe$
(with $|g| = 5$).  Now $T_4C$ is a tensor algebra on $e$, $f$, and
$g$, with certain differentials.  One next looks at the homology in
dimension five, and continues.  Clearly this process gets very
cumbersome as one goes higher and higher in the resolution.  (Note: A
complete Koszul-type description of the resolution in this case has
recently been given by Bill Kronholm).
\end{example}

\subsection{Homotopy classes}
\label{sec:homotopy}

Suppose that $C$ and $D$ are $k$-dgas, and that $C$ is cofibrant.  To
compute maps in the homotopy category $\ho(\kDGA)(C,D)$, one may
either use a cylinder object for $C$ or a path object for
$D$~\cite[1.2.4, 1.2.10]{hovey}.  In the case of general dgas, both
are somewhat complicated.  For the very special situations that arise
in this paper, however, there is a simple method using path objects.

Assume $C_i=0$ for $i<0$, $C_0$ is generated by $1$ as a $k$-algebra,
and that $D=k\oplus M$ for some $M\in \Ch_k$.  Here $D$ is the dga
obtained by adjoining $M$ to $k$ as a square-zero ideal.  Let $I$
denote the chain complex $k \ra k^2$ concentrated in degrees $1$ and
$0$, where $d(a)=(a,-a)$.  Recall that $\Hom_k(I,M)$ is a path object
for $M$ in $\Ch_k$, and let $PD$ denote the square-zero extension
$k\oplus \Hom_k(I,M)$.  This is readily seen to be a path object for
$D$ in $\kDGA$.  It is not a {\it good\/} path object, however, as the
map $PD \ra D\times D$ is not surjective in degree $0$.  Despite this
fact, it is still true that $\ho(\kDGA)(C,D)$ is the coequalizer of
$\kDGA(C,PD)\dbra \kDGA(C,D)$; this is a simple argument using that all
maps of dgas with domain $C$ coincide on $C_0$.

\subsection{Ring spectra}
There are, of course, different settings in which one can study ring
spectra.  We will work with the category of symmetric spectra
$\Sp^\Sigma$ of \cite{hss} with its symmetric monoidal smash product
$\Smash$ and unit $S$.  A \dfn{ring spectrum} is just an
$S$-algebra---that is to say, it is a spectrum $R$ together with a
unit $S\ra R$ and an associative and unital pairing $R\Smash R \ra R$.
The category $\RingSpectra$ has a model structure in which a map is a
weak equivalence or fibration precisely if it is so when regarded as a
map of underlying spectra.  See \cite[5.4.3]{hss}. 

The forgetful functor $\RingSpectra \ra \Sp^\Sigma$ has a left adjoint
$T$.  So for any spectrum $E$ there is a `free ring spectrum built
from $E$', denoted $T(E)$.  The cofibrations in $\RingSpectra$ are
generated (via retracts, cobase-change, and transfinite composition)
from those of the form $T(A) \ra T(B)$ where $A\ra B$ ranges over the
generating cofibrations of $\Spectra$.  Recall that for $\Spectra$
these generators are just the maps $\Sigma^\infty(\bdd{n})\ra
\Sigma^\infty(\del{n})$ and their desuspensions.  The situation
therefore exactly parallels that of dgas.

\subsection{Eilenberg-Mac\,Lane ring spectra}
\label{se:EM}
Any dga $A$ gives rise to a ring spectrum $HA$ called the
\dfn{Eilenberg-Mac\,Lane ring spectrum associated to $A$}.  This can be
constructed functorially, and has the property that if $A\ra B$
is a quasi-isomorphism then $HA\ra HB$ is a weak equivalence.
It is also the case that $H$ preserves homotopy limits.

Unfortunately, giving a precise construction of $H(\blank)$ seems to
require a morass of machinery.  This is accomplished in \cite{S}.
We will give a brief summary, but the reader should note that these
details can largely be ignored for the rest of the paper.

Let $\ch_+$ be the category of non-negatively graded chain complexes.
On this category one can form symmetric spectra based on the object
$\bZ[1]$, as in \cite{stab}; call the corresponding category
$\Sp^\Sigma(\ch_+)$.  Similarly one can form symmetric spectra
based on simplicial abelian groups and the object $\widetilde{\bZ}S^1$
(the reduced free abelian group on $\Delta^1/\bdd{1}$);
call this category $\Sp^\Sigma(s\Ab)$.

There are two Quillen equivalences (with left
adjoints written on top)
\begin{myequation}
\label{eq:adjoint}
\xymatrix{
\Sp^\Sigma(\ch_+) \ar@<0.5ex>[r]^-D & \Ch_{\bZ} \ar@<0.5ex>[l]^-R
& \text{and} &
   \Sp^\Sigma(\ch_+) \ar@<0.5ex>[r]^-L &\Sp^\Sigma(s\Ab)
\ar@<0.5ex>[l]^-{\namma}
}
\end{myequation}
in which $(D, R)$ is strong monoidal and $(L, \namma)$ is weak
monoidal (the functor $\namma$ is denoted $\phi^*N$ in \cite{S}).  See
\cite[3.6]{SS3} 
 for the terms `strong monoidal' and `weak
monoidal'.  By the work in \cite{S} the above functors induce Quillen
equivalences between the corresponding model categories of ring
objects (or monoids).  Thus we have adjoint pairs
\[ D\colon \Ring(\Sp^\Sigma(\ch_+)) \adjoint \Ring(\Ch_{\bZ})\colon R \qquad
\]
\[
L^{mon}\colon \Ring(\Sp^\Sigma(\ch_+)) \adjoint \Ring(\Sp^\Sigma(s\Ab))
\colon \namma.
\]
In the first case the functors $D$ and $R$ are just the restriction of
those in (\ref{eq:adjoint}), as these were strong monoidal.  But in
the second case only the right adjoint is restricted from
(\ref{eq:adjoint}); the left adjoint is more complicated.
See \cite[3.3]{SS3} 
 for a complete description of $L^{mon}$.

Finally, the Quillen pair $F\colon \sSet \adjoint s\Ab\colon U$
(where $U$ is the forgetful functor) induces a strong monoidal
Quillen pair $F\colon \Sp^\Sigma \adjoint \Sp^\Sigma(s\Ab)\colon U$.
Therefore, by \cite[3.2]{SS3}, 
 there is a Quillen pair
\[F\colon \RingSpectra \adjoint
\Ring[\Sp^\Sigma(s\Ab)]\colon U\]

Let $\und{R}$, $\und{D}$, etc. denote the derived functors of $R$
and $D$---that is, the induced functors on homotopy categories.  If
$A$ is a dga, we then define
\[ HA= \und{U} \bigl [ \und{L}^{mon} (\und{R} A) \bigr ].
\]
This is a ring spectrum with the desired properties.  To see that $H$
preserves homotopy limits, for instance, note that both $L^{mon}$ and
$R$ have this property because they are functors in a Quillen
equivalence; likewise $U$ has this property because it is a right
Quillen functor.

It is also convenient to relate the above constructions to the
categories of $H\Z$-algebras and modules.
Since the unit in $\Sp^\Sigma(s\Ab)$ forgets via $U$ to the symmetric spectrum
$H\Z$, the Quillen pair $(F,U)$ factors through a strong monoidal
Quillen equivalence $Z\colon H\Z\lmod \adjoint \Sp^\Sigma(\sab)\colon
U'$ by~\cite[2.5]{S}.  Applying \cite[3.2]{SS3} 
again yields
a Quillen equivalence
\[ Z\colon H\Z\Algebra \adjoint \Ring[\Sp^\Sigma(s\Ab)]\colon U',\]
and the functor $U\colon \Ring[\Sp^\Sigma(s\Ab)] \ra S\Alg$ is the
composite of $U'$ with the restriction $H\Z\Alg \ra S\Alg$.
This shows that $HA$ is naturally an $H\Z$-algebra.

\section{Postnikov sections and $k$-invariants for dgas}
\label{se:Postnikov-dga}
Fix a commutative ring $k$.  In this section we work in the model
category $\kDGA$, and describe a process for understanding the
quasi-isomorphism type of a dga.  This involves building the dga from
the ground up, one degree at a time, by looking at its Postnikov
sections.  The difference between successive Postnikov sections is
measured by a homotopical extension, usually called a $k$-invariant
(unfortunately the `$k$' in `$k$-invariant' has no relation to the
commutative ring $k$!)  This $k$-invariant naturally lives in a
certain Hochschild cohomology group.  In Example~\ref{ex:exterior} we
use these tools to classify all dgas over $\Z$ whose homology is
$\Lambda_{\F_p}(g_n)$.

\subsection{Postnikov sections}

\begin{defn}\label{3.2}
If $C$ is a dga and $n\geq 0$, an \mdfn{$n$th Postnikov section} of
$C$ is a dga $X$ together with a map $C\ra X$ such that
\begin{enumerate}[(i)]
\item $H_i(X)=0$ for $i>n$, and
\item $H_i(C) \ra H_i(X)$ is an isomorphism for $i\leq n$.
\end{enumerate}
\end{defn}

Fix an $n\geq 0$.  We are able to construct Postnikov sections only
when $C$ is a \dfn{non-negative} dga---meaning that $C_i=0$ when
$i<0$.  A functorial $n$th Postnikov section can be obtained by
setting
\[
[\PP_n C]_i=
\begin{cases} C_i & \text{if $i < n$}, \\
C_n/\im(C_{n+1}) & \text{if $i=n$, and} \\
0 & \text{if $i>n$.}
\end{cases}
\]
This construction is sometimes inconvenient in that the map $C\ra
\PP_n C$ is not a cofibration.  To remedy this, there is an
alternative construction using the small object argument.  Given any
cycle $z$ of degree $n+1$, there
is a unique map $k(S^{n+1}) \ra C$ sending the generator in degree
$n+1$ to $z$.  We can construct the pushout of $k(D^{n+2}) \la
k(S^{n+1}) \ra C$, which has the effect of killing $z$ in homology.
Let $L_{n+1}C$ be the pushout
\[ \coprod k(D^{n+2}) \lla \coprod k(S^{n+1}) \lra C \]
where the coproduct runs over all cycles in degree $n+1$.
Define
$P_nC$ to be the colimit of the sequence
\[ C \ra L_{n+1}C \ra L_{n+2}L_{n+1}C \ra L_{n+3}L_{n+2}L_{n+1}C \ra
\cdots \]
It is simple to check that this is indeed another functorial
$n$th Postnikov section,
and $C\ra P_nC$ is a cofibration.  Note as well that there is a
natural projection $P_nC \ra \PP_nC$, which is a quasi-isomorphism.

\begin{prop}\label{prop-can-postnikov}
For any non-negative dga $C$ and any $n$th Postnikov section $X$, there is
a quasi-isomorphism $P_nC \to X$.
\end{prop}

\begin{proof}
The quasi-isomorphism can be constructed directly using the description
of the functors $L_i$.
\end{proof}

It follows from the above proposition that any two $n$th Postnikov
sections for $C$ are quasi-isomorphic.

Note that there are canonical maps $P_{n+1}C \ra P_nC$ compatible with
the co-augmentations $C\ra P_{i}C$.
The sequence $\cdots \ra P_{n+1}C \ra P_nC \ra \cdots \ra P_0C$ is
called the \dfn{Postnikov tower} for $C$.

\subsection{k-invariants}
Let $C$ be a dga and let $M$ be a $C$-bimodule---i.e., a $(C\tens_k
C^{op})$-module.  By $C\Wedge M$ we mean the square-zero extension of
$C$ by $M$; it is the dga whose underlying chain complex is $C\oplus
M$ and whose algebra structure is the obvious one induced from the
bimodule structure (and where $m\cdot m'=0$ for any $m,m'\in M$).

Assume $C$ is non-negative.  Then there is a natural map $P_0(C) \ra
H_0(C)$, and it is a quasi-isomorphism (here $H_0(C)$ is regarded as a
dga concentrated in degree zero).  Since $H_{n+1}(C)$ is a
bimodule over $H_0(C)$, it becomes a bimodule over each $P_nC$ by
restriction through $P_nC \ra P_0C \ra H_0(C)$.  So we can look at the
square-zero extension $P_nC \Wedge \Sigma^{n+2}H_{n+1}(C)$

There is a canonical map of dgas $\gamma\colon P_{n}C \ra \PP_{n}C \Wedge
\Sigma^{n+2} [H_{n+1}C]$ defined by letting it be the identity in
dimensions smaller than $n$, the natural projection in dimension $n$,
and the zero map in dimension $n+1$ and all dimensions larger than
$n+2$.  In dimension $n+2$ it can be described as follows.  The map is
zero on $C_{n+2}$, and if $x\in [P_nC]_{n+2}$ was adjoined to kill the
cycle $z\in C_{n+1}$, then $x$ is mapped to the class of $z$ in
$\Sigma^{n+2}[H_{n+1}(C)]$.  A little checking shows $\gamma$ is a
well-defined map of dgas.

Let $\Ho(\kDGA_{/{\PP_nC}})$ denote the homotopy category of $k$-dgas
augmented over $\PP_nC$, and let $\alpha_n \in
\Ho(\kDGA_{/\PP_nC})(P_nC,\PP_nC \Wedge \Sigma^{n+2}[H_{n+1}C])$ be the
image of the map $\gamma$.
Then $\alpha_n$ is called the \mdfn{$n$th
$k$-invariant} of $C$.  One can check that it depends only on the
quasi-isomorphism type of $P_{n+1}C$.
Moreover, the homotopy type of
$P_{n+1}C$ can be recovered from $\alpha_n$.  This is shown by the
following result, since of course the homotopy fiber of $\gamma$ only
depends on its homotopy class:

\begin{prop}
$P_{n+1}C \ra P_n C \llra{\gamma} \PP_{n}C
\Wedge \Sigma^{n+2}[H_{n+1}C]$ is a homotopy fiber sequence in
$\Ho(\kDGA_{/\PP_nC})$.
\end{prop}

\begin{proof}
The result may be rephrased as saying that $P_{n+1}C$ is weakly
equivalent to the homotopy pullback (in the category of dgas) of the
diagram
\[ \xymatrix{
P_nC \ar[r] & \PP_nC\Wedge \Sigma^{n+2}[H_{n+1}C] & \PP_nC \ar@{
>->}[l]
}
\]
where the right map is the obvious inclusion.  Note also that the left
map is a fibration, by construction, and so the pullback and homotopy
pullback are equivalent.  One readily sees that there is a natural map
from $P_{n+1}C$ to this pullback, and that it is a quasi-isomorphism.
\end{proof}

\begin{remark}
As $\kDGA$ is a right proper model category, it follows that for any
 $X\we Y$ the induced Quillen map $\kDGA_{/X} \ra
\kDGA_{/Y}$ is a Quillen equivalence.  The quasi-isomorphism $P_n C\ra
\PP_n C$ therefore allows us to identify the set
$\Ho(\kDGA_{/\PP_nC})(P_nC,\PP_nC \Wedge \Sigma^{n+2}[H_{n+1}C])$ with
\[ \Ho(\kDGA_{/P_nC})(P_nC,P_nC \Wedge \Sigma^{n+2}[H_{n+1}C]).
\]
\end{remark}

\subsection{Classifying extensions}
\label{se:classifyext}

\begin{defn}
Let $C$ be a non-negative $k$-dga such that $H_i(C)=0$ for $i> n$.  A
\mdfn{Postnikov $(n+1)$-extension} of $C$ is a $k$-dga $X$ such that
$P_{n+1}X\he X$ together with a map of $k$-dgas $f\colon X \ra C$
which yields an isomorphism on $H_i(\blank)$ for $i\leq n$.  A
map $(X,f) \ra (Y,g)$ between Postnikov $(n+1)$-extensions is
defined to be a quasi-isomorphism $X\ra Y$ compatible with $f$ and
$g$.
\end{defn}

\begin{prop}
\label{pr:k-inv}
Assume $C$ is non-negative and that $P_nC \he C$.  Fix an
$H_0(C)$-bimodule $M$.  Consider the category whose objects consist of
Postnikov $(n+1)$-extensions $(X,f)$ together with an isomorphism of
$H_0(C)$-bimodules $\theta\colon H_{n+1}(X) \ra M$.  A map from
$(X,f,\theta)$ to $(Y,g,\sigma)$ is a quasi-isomorphism $X\ra Y$
compatible with the other data.  Then the connected components of this
category are in bijective correspondence with the set of homotopy
classes $\Ho(\kDGA_{/C})(C,C \Wedge \Sigma^{n+2}M)$.
\end{prop}

In the context of the above result, the general problem one would like
to be able to solve is to classify all Postnikov $(n+1)$-extensions
$X$ of $C$ having $H_{n+1}X \iso M$.  It is not quite true that the
set of quasi-isomorphism types of all such $X$ is in bijective
correspondence with $\Ho(\kDGA_{/C})(C,C \Wedge \Sigma^{n+2}M)$.
Different $k$-invariants can nevertheless lead to quasi-isomorphic
dgas $X$.  To see this, note that if $h\colon M \ra M$ is an
automorphism of $H_0C$-bimodules then any $k$-invariant can be
`twisted' by this automorphism.  The homotopy fibers of the original
and twisted $k$-invariants are quasi-isomorphic, however---the
only thing that is different about them is the prescribed isomorphism
between their $H_{n+1}$ and $M$.  This is why such an isomorphism must
be built into the category appearing in the proposition.

\begin{proof}[Proof of Proposition~\ref{pr:k-inv}]
Let $\cA$ be the category described in the statement of the
proposition, and let $T=\Ho(\kDGA_{/C})(C,C\Wedge
\Sigma^{n+2}M)$.

Suppose $X\ra C$ is a Postnikov $(n+1)$-extension of $C$, and $\theta
\colon H_{n+1}X \ra M$ is an isomorphism of $H_0(C)$-bimodules.  As
above, construct the map $\gamma\colon P_n X \ra \PP_nX \Wedge
\Sigma^{n+2}[H_{n+1}X]$.  The map $X\ra C$ induces quasi-isomorphisms
$P_nX \ra P_n C$ and $\PP_n X \ra \PP_n C$.  Note that the maps $C\ra
P_n C$ and $C\ra \PP_n C$ are quasi-isomorphisms as well.  These,
together with $\theta$, allow us to identify $\gamma$ with a map in
the homotopy category $C \ra C \Wedge \Sigma^{n+2}M$.  One checks that
this gives a well-defined map $\pi_0\cA \ra T$.

Now suppose $\alpha\in T$.  Let $\tilde{C} \ra C$ be a
cofibrant-replacement, and let $\tilde{C} \ra C\Wedge \Sigma^{n+2}M$
be any map representing $\alpha$.  Let $X$ be the homotopy pullback of
$\tilde{C} \ra C \Wedge \Sigma^{n+2}M \la C$ where the second map is
the obvious inclusion.  The composition $X\ra \tilde{C} \ra C$ makes
$X$ into a Postnikov $(n+1)$-extension of $C$, and the long exact
sequence for the homology of a homotopy pullback gives us an
isomorphism $\theta\colon H_{n+1}X \ra M$.  One checks that this
defines a map $T \ra \pi_0\cA$.  With some trouble one can verify that
this is inverse to the previous map $\pi_0\cA \ra T$.
(Remark~\ref{re:post-paper} below suggests a better, and more
complete, proof).
\end{proof}

The following corollary of the above proposition is immediate:

\begin{cor}
\label{co:classify}
Let $C$ and $M$ be as in the above proposition.
Let $G$ be the group of $H_0(C)$-bimodule automorphisms of $M$.  Let $S$
be the quotient set of $G$
 acting on $\Ho(\kDGA_{/C})(C,C \Wedge \Sigma^{n+2}M)$.
Consider the category of Postnikov $(n+1)$-extensions $X$
of $C$ which satisfy $H_{n+1}X \iso M$ as $H_0(C)$-bimodules (but
where no prescribed choice of isomorphism is given).
Then $S$ is in bijective correspondence with the connected
components of this category.
\end{cor}

\begin{proof}
Let $\cA$ be the category described in Proposition~\ref{pr:k-inv} and
let $\cB$ be the category described in the statement of the
corollary.  There is clearly a surjective map $\pi_0(\cA)/G \ra
\pi_0(\cB)$.  Injectivity is a very simple exercise.
\end{proof}

\begin{remark}
\label{re:post-paper}
A more complete proof of the above proposition and corollary can be
obtained by following the methods of \cite{post}.  That paper takes
place
in the setting of ring spectra, but all the arguments adapt verbatim.
Alternatively,
using~\cite{S} one can consider $k$-dgas as $Hk$-algebra spectra---so
from that perspective the above results are actually special cases of
those
from~\cite[8.1]{post}.  
\end{remark}

\begin{remark}\label{rem-conn}
The above material can be applied to any connective dga $C$ (that is,
one where $H_k(C)=0$ for $k<0$).  Such a dga is always
quasi-isomorphic to a non-negative dga $QC$.
One gets $k$-invariants in the set
\[
\Ho(\kDGA_{/P_n(QC)})\bigl (P_n(QC),P_n(QC)\Wedge
\Sigma^{n+2}H_{n+1}(QC) \bigr ).
\]
\end{remark}

\begin{example}
\label{ex:H=ext}
Let $C$ be a dga over $\Z$ with $H_*(C)$ equal to an exterior algebra
over $\Fp$ on a generator in degree $2$.  What are the possibilities
for $C$?  We know that $P_1(C)\he \Fp$ and $P_2(C)\he C$.  We therefore
have to analyze the single homotopy fiber sequence $P_2C \ra \Fp \ra \Fp
\Wedge \Sigma^3 \Fp$.  What are the possibilities
for the $k$-invariant in this sequence?  One has to remember here that
$\Fp$ is not cofibrant as a dga over $\Z$, and so one must work with
a cofibrant replacement.

The first few degrees of a cofibrant replacement for $\Fp$ look like
\[ \xymatrix{ \cdots \ar[r] & \Z e^3 \oplus \Z f \ar[r]^-{(p,1)} &\Z
e^2 \ar[r]^0 &\Z e
\ar[r]^{p} & \Z.1
}
\]
These symbols mean $d(e)=p$ (which implies $d(e^n)=pe^{n-1}$ when $n$
is odd and $d(e^n)=0$ when $n$ is even) and $d(f)=e^2$.  We are
interested in maps from this dga to the simpler dga
\[ \xymatrix{ 0 \ar[r]& (\Z/p).g\ar[r] & 0\ar[r] & 0\ar[r] & (\Z/p).1 }
\]
(an exterior algebra with a class in degree $3$, and zero
differential).  The possibilities for such maps are clear: $e$ must be
sent to $0$, and $f$ can be sent to a (possibly zero) multiple of $g$.
One finds that $\Ho(\ZDGA_{/\Fp})(\Fp,\Fp\Wedge \Sigma^3 \Fp) \iso
\Z/p$; this requires an analysis of homotopies, but using
(\ref{sec:homotopy}) one readily sees that all homotopies are
constant.  Now, the group of automorphisms of $\Fp$ as an $\Fp$-module
is just $(\Z/p)^*$, and there are precisely two orbits of $\Z/p$ under
this group action.  By Corollary~\ref{co:classify}, we see that there
are precisely two quasi-isomorphism types of $\Z$-dgas  having
homology algebra $\Lambda_{\Fp}(g_2)$.

We can find these two dgas explicitly by constructing the appropriate
homotopy pullbacks.  When the $k$-invariant has $f \mapsto 0$ one
finds that $C$ is just an exterior algebra with zero differential (the
$k$-invariant $\Fp \ra \Fp \Wedge \Sigma^3 \Fp$ is just the obvious
inclusion).  When the $k$-invariant has $f\mapsto g$ one has that $C$
is quasi-isomorphic to the dga $\Z[e; de=p]/(e^4)$.
\end{example}

\subsection{Hochschild cohomology}
We want to extend Example~\ref{ex:H=ext} and classify all
dgas whose homology algebra is $\Lambda_{\Fp}(g_n)$.  A problem
arises, in that one has to compute a cofibrant-replacement of $\Fp$
(as a $\Z$-dga) up to dimension $n+1$.  Such a cofibrant-replacement
becomes very large in high dimensions.  Hochschild cohomology gives a
way around this issue, which we now recall.

One possible reference for the material in this section is
\cite[Section 2]{L}.  Lazarev works in the context of ring spectra,
but all his results and proofs work exactly the same for dgas.
Looked at differently, $k$-dgas are the same as ring spectra which are
$Hk$-algebras by~\cite{S}---so the results below are just special cases of
Lazarev's, where the ground ring is $Hk$.

If $C$ is a $k$-dga, let $\Omega_{C/k}$ denote the homotopy fiber of
the multiplication map $C \tens^L_{k} C \ra C$.  From now on we will
just write $\tens_k$ instead of $\tens^L_k$; but it is important to
never forget that all tensors are now {\it derived\/} tensors.  We
also write $\Hom$ rather than $\RHom$.  If $M$ is a $C\tens_k
C^{op}$-module, there is an induced homotopy fiber sequence
\[ \Fhom_{C\tens_k C^{op}}(C,M) \ra \Fhom_{C\tens_k C^{op}}(C\tens_k C^{op},M) \ra
\Fhom_{C\tens_k C^{op}}(\Omega_{C/k},M). \]
The term in the middle may be canonically identified with $M$.
One typically makes the following definitions:
 \[ \sDer_k(C,M)=\Fhom_{C\tens_k C^{op}}(\Omega_{C/k},M), \qquad
\sHH_k(C,M)=\Fhom_{C\tens_k C^{op}}(C,M) \]
\[ \Der^n_k(C,M)=H_{-n}[\sDer_k(C,M)], \qquad\text{and}\qquad
   \HH^n_k(C,M)=H_{-n}[\sHH_k(C,M)].
\]
Sometimes we will omit the $k$ subscripts when the ground ring is
understood.

The homotopy fiber sequence $\sHH(C,M) \ra M \ra
\sDer(C,M)$ gives a long exact sequence of the form
\[ \cdots \ra \HH^n(C,M) \ra H_{-n}(M) \ra \Der^n(C,M) \ra
HH^{n+1}(C,M) \ra \cdots
\]
We will mostly be interested in applying this when $M$ is concentrated
in a single dimension $r$, in which case $\Der^{*}(C,M)\iso
\HH^{*+1}(C,M)$ for $*\notin \{-r,-r-1\}$.

Finally, in order to connect all this with the classification of dgas,
one can prove that there is an isomorphism
\[\label{der=ho} \Der^n(C,M) \iso
\Ho(\kDGA_{/C})(C,C\Wedge \Sigma^n M).
\]
The proof of this isomorphism in \cite{L} seems to contain gaps; we are
very grateful to Mike Mandell for showing us a complete proof \cite{M}.

\begin{example}
\label{ex:exterior}
We'll use the above machinery to determine all dgas $C$ over $\Z$ such
that $H_*(C)\iso \Lambda_{\Fp}(g_n)$ (an exterior algebra on a class
of degree $n$).  Such a dga has $P_{n-1}(C)\he \Fp$ and $P_nC \he C$,
so we have the homotopy fiber sequence $C \ra \Fp \ra \Fp \Wedge
\Sigma^{n+1} \Fp$.  We need to understand the possibilities for the
second map.

The above observations give us isomorphisms
\[ \Ho(\Z-\DGA_{/\Fp})(\Fp,\Fp \Wedge \Sigma^{n+1} \Fp) \iso
\Der^{n+1}_\Z(\Fp,\Fp) \iso \HH^{n+2}_\Z(\Fp,\Fp) \]
where for the last isomorphism we need $n \notin \{-1,-2\}$.
Recall furthermore that $\HH^{*}(\Fp,\Fp) \iso H_{-*}[\Fhom_{\Fp \tens
\Fp}(\Fp,\Fp)]$, and remember that $\Fp \tens \Fp$ really means $\Fp
\tens^L_{\Z} \Fp$ here.

The dga $\Fp \tens^L \Fp$ is $\Lambda=\Lambda_{\Fp}(e_1)$ (an
exterior algebra with zero differential).  Of course
$\Ext_{\Lambda}(\Fp,\Fp)$ has homology algebra equal to a polynomial
algebra on a class of degree $-2$ (or $+2$ if  cohomological
grading is used).  That is, $\HH^*(\Fp,\Fp) \iso
\Fp[\sigma_2]$. 

The conclusion is that when $n \geq 1$ is odd, there is only one
homotopy class in $\Ho(\Z-\DGA_{/\Fp})(\Fp,\Fp \Wedge \Sigma^{n+1}
\Fp)$ (the trivial one), and hence only one quasi-isomorphism type for
dgas whose homology is $\Lambda_{\Fp}(g_n)$ (given by this graded
algebra with zero differential).  
When $n\geq 2$ is even we have
\[ \Ho(\Z-\DGA_{/\Fp})(\Fp,\Fp \Wedge \Sigma^{n+1} \Fp)\iso \Z/p,
\]
and
quotienting by $\Aut(\Fp)=(\Fp)^*$ gives exactly two orbits.  So in
this case there are exactly two such dgas: the trivial square-zero
extension (exterior algebra) and an `exotic' one.

For example, when $n=2$ the non-trivial dga is
$\Z[e;de=p]/(e^3,pe^2)$ with $|e| = 1$.  When $n=4$ it is 
$\Z\langle e,f;de=p,
df=e^2 \rangle/(e^4, e^2f, efe, fe^2, ,fef, f^2, p(ef + fe))$
with $|e| = 1, |f| = 3$. 
\end{example}

\begin{example}
Suppose that $F$ is a field, and we want to classify all $F$-dgas
whose homology is $\Lambda_F(g_n)$.  Everything proceeds as above, and
we find ourselves needing to compute $\HH^*_F(F,F)$.  But this is
trivial except when $*=0$, and so there is only one quasi-isomorphism
type for the dgas in question---namely, the trivial one.  The previous
example is more complicated because the ground ring is $\Z$.
\end{example}


\section{$k$-invariants for ring spectra}
\label{se:Postnikov-ringsp}
The material developed for dgas in the last section is
developed for ring spectra in~\cite{post}.
One has Postnikov
towers of ring spectra, and the $k$-invariants measuring the
extensions at each level now live in groups called {\it topological
Hochschild cohomology\/}.  We use these tools to classify ring spectra
whose homotopy algebra is $\Lambda_{\F_p}(g_n)$.
Another reference for some of
this background material is~\cite[Sections 2 and 8.1]{L}.

\smallskip

Fix a commutative ring spectrum $R$.

\begin{defn}
If $T$ is an $R$-algebra and $n\geq 0$, an \mdfn{$n$th Postnikov
section} of $T$ is an $R$-algebra $P_nT$ together with a map of
$R$-algebras $T\ra P_n T$ such that
\begin{enumerate}[(i)]
\item $\pi_i(P_nT)=0$ for $i>n$, and
\item $\pi_i(T) \ra \pi_i(P_nT)$ is an isomorphism for $i\leq n$.
\end{enumerate}
\end{defn}

Postnikov sections can be constructed for connective $R$-algebras just
as they were for dgas.  See \cite[Section 2]{post} for details.
A connective $R$-algebra $T$ has $k$-invariants lying in the set of homotopy
classes $\Ho(R\Alg_{/(P_nT)})(P_nT, P_nT \Wedge \Sigma^{n+2}
H(\pi_{n+1}T))$, giving homotopy fiber sequences
\[ P_{n+1}T \ra P_n T \ra P_nT\Wedge \Sigma^{n+2}H(\pi_{n+1}T).
\]
This is where things diverge somewhat from what we did for dgas.  For
dgas we were able to produce the $k$-invariants in terms of very
explicit formulas, defined on the elements of the dga.  One cannot use
this construction for ring spectra.  Instead, one has to use a more
categorical approach.  This is developed in~\cite{post}
and the analog of Corollary~\ref{co:classify} is proven in~\cite[8.1]{post}.

If $M$ is a $T\Smash_R T^{op}$-module, one defines $\sTDer_R(T,M)$ and
$\sTHH_R(T,M)$ just as in the previous section.

\begin{example}
Let us investigate all ring spectra (i.e., $S$-algebras) $T$ whose
homotopy algebra is $\pi_*T\iso \Lambda_{\Fp}(g_n)$ ($n\geq 1$).  Just
as for dgas, we have the single homotopy fiber sequence $T \ra H\Fp
\ra H\Fp \Wedge \Sigma^{n+1} H\Fp$ in $\ho(S\Alg_{/H\F_p})$.  The
possibilities for the $k$-invariant are contained in the set
$\TDer^{n+1}(H\Fp,H\Fp)\iso\THH^{n+2}(H\Fp,H\Fp)$.  These
$\THH$-groups have been calculated by B\"okstedt \cite{Bok}, but see
\cite[5.2]{HM} 
 for a published summary (those references deal
with $THH$ {\it homology\/}, but one can get the cohomology groups by
dualization).  We have $\THH^*(\Fp,\Fp)\iso \Gamma[\alpha_2]$, a
divided polynomial algebra on a class of degree $2$.
As another source for this computation, including the ring structure,
we refer to \cite[7.3]{FLS}. 

It follows from the computation that when $n\geq 1$ is odd there is
only one ring spectrum with the given homotopy algebra, namely the
Eilenberg-Mac\,Lane spectrum
$H(\Lambda_{\Fp}(g_n))$.  When $n\geq 2$ is even, there are exactly
two such ring spectra.
\end{example}

\begin{remark}
In the above example, we'd like to call special attention to the case
where $n=2p-2$ for $p$ a prime.  In this case we can say precisely
what the two homotopy types of ring spectra are.  One of them, of
course, is the trivial example $H(\Lambda_{\Fp}(g_{2p-2}))$.  The
other is the first nontrivial Postnikov section of connective Morava
$K$-theory, $P_{2p-2}k(1)$ (see \cite{A} or \cite{G} for the ring
structure on $k(1)$ when $p=2$).  These two ring spectra are obviously
not weakly equivalent, since their underlying spectra are
not even weakly equivalent (the latter is not an Eilenberg-Mac\,Lane
spectrum).  So they must represent the two homotopy types.
\end{remark}

\subsection{Comparing $\HH$ and $\THH$}
Suppose $Q\ra R$ is a map of commutative ring spectra, and $T$ is an
$R$-algebra.  There is a natural map $T\Smash_Q T^{op} \ra T \Smash_R
T^{op}$, and as a result a natural map
\[\sTHH_R(T,M)=\Fhom_{T\Smash_R T^{op}}(T,M) \ra \Fhom_{T\Smash_Q T^{op}}(T,M)
=\sTHH_Q(T,M).
\]
(Note that the smash means `derived smash' and the hom means
`derived hom', as will always be the case in this paper).
In particular we may apply this when $Q\ra R$ is the map $S \ra H\Z$.
If $T$ is an $H\Z$-algebra we obtain a map $\sTHH_{H\Z}(T,M) \ra
\sTHH_S(T,M)$.  By the equivalence of $H\Z$-algebras with dgas,
$\sTHH_{H\Z}$ is just another name for $\sHH_{\Z}$.  So if $T$ is a
dga and $M$ is a $T$-bimodule, we are saying that there are natural
maps
\[ \HH^*_\Z(T,M) \ra \THH^*_S(HT,HM).
\]
When $M=T$ both the domain and codomain have ring structures, and this
is a ring map.

When $M$ is concentrated entirely in degree $0$ (as in all our
application), one can show that for $n\geq 1$ the following square
commutes:
\[
\xymatrix{ \ho(H\Z\Alg_{/HT})(HT,HT\Wedge \Sigma^{n+1} HM) \ar[r]^-\iso\ar[d]
& \THH^{n+2}_{H\Z}(T,M) \ar[d] \\
\ho(S\Alg_{/HT})(HT,HT\Wedge \Sigma^{n+1} HM) \ar[r]^-\iso
& \THH^{n+2}_S(HT,HM).
}
\]
Then, using the identification of $H\Z$-algebras with dgas,
the top horizontal map can be identified with
\[
\ho(\Z\!-\!\DGA_{/T})(T,T\Wedge \Sigma^{n+1} M) \to
\HH^{n+2}_\Z(T,M).
\]

\smallskip

\section{Examples of topological equivalence}
\label{se:examples}

In this section we present several examples of dgas which are
topologically equivalent but not quasi-isomorphic.  Since the homotopy
theory of dgas is Quillen equivalent to that of $H\Z$-algebras, this
is the same as giving two non-equivalent $H\Z$-algebra structures on
the same underlying ring spectrum.  Our examples are:
\begin{enumerate}[(a)]
\item The dgas $\Z[e_1;de_1=2]/(e_1^4)$ and $\Lambda_{\F_2}(g_2)$.
\item The two distinct quasi-isomorphism types of dgas
whose homology is the exterior algebra $\Lambda_{\F_p}(g_{2p-2})$, provided by
Example~\ref{ex:exterior} (here $p$ is a fixed prime).  
This gives the example from (a) when $p=2$.
\item The ring spectrum $H\Z \Smash H\Z/2$, with the two structures of
$H\Z$-algebra provided by the two maps $H\Z = H\Z\Smash S \ra
H\Z\Smash H\Z/2$ and $H\Z=S\Smash H\Z \ra H\Z \Smash H\Z/2$.
\item The ring spectrum $H\Z \Smash_{bo} H\Z/2$, with the two structures of
$H\Z$-algebra coming from the left and the right as in (c).
\end{enumerate}

For parts (a) and (b), we must prove that these dgas are topologically
equivalent---we do this by using the comparison map from Hochschild
cohomology to topological Hochschild cohomology.  For parts (c) and
(d), we must prove that the associated dgas are not
quasi-isomorphic---we do this by calculating the derived tensor with
$\Z/2$, and finding that we get different homology rings.

The examples in (a), (c), and (d) are related.  Specifically, the two
dgas in (a) are the second Postnikov sections of the dgas in (d) (or
in (c)).

\subsection{Examples using $\HH$ and $\THH$}
\label{ex:tpwe}
In this section we will mainly be working with
$\Z$-dgas and with $S$-algebras.  The symbols $\HH^*$ and $\THH^*$ will
always indicate $\HH^*_{\Z}$ and $\THH^*_{S}$, unless otherwise noted.

We begin with the following simple result:

\begin{prop} \label{prop-HPnC}
Let $C$ be a dga, and let $P_nC$ be an $n$th Postnikov section.  Then
$HC \ra H(P_nC)$ is an $n$th Postnikov section for the ring spectrum
$HC$.
\end{prop}

\begin{proof}
Immediate.
\end{proof}

Suppose $C$ is a non-negative $\Z$-dga with $P_nC\he C$, and let $X$ be
a Postnikov $(n+1)$-extension of $C$.  Write $M=H_{n+1}X$.  We have a
homotopy fiber sequence $X \ra C \ra C \Wedge \Sigma^{n+2}M$ in
$\Ho(\Z-DGA_{/C})$.  Since $H(\blank)$ preserves homotopy limits,
applying it yields a homotopy fiber sequence $HX \ra HC \ra HC \Wedge
\Sigma^{n+2}HM$ in $\Ho(\RSp{/HC})$.  So $HX$ is a Postnikov
$(n+1)$-extension of $HC$.  The $k$-invariant for $X$ lies in
$\HH^{n+3}(C,M)$, whereas the $k$-invariant for $HX$ lies in
$\THH^{n+3}(HC,HM)$.  The latter is the image of the former under the
natural map $\HH^{n+3}(C,M) \ra \THH^{n+3}(HC,HM)$.

We can now give our first example of two dgas which are topologically
equivalent but not quasi-isomorphic.  The example will be based on our
knowledge of the map $\HH^*(\Fp,\Fp) \ra \THH^*(\Fp,\Fp)$.  The domain
is $\Fp[\sigma_2]$ (as calculated in \ref{ex:exterior}), whereas the codomain
is the divided power algebra $\Gamma_{\Fp}[\alpha_2]$
(cf. \cite[7.3]{FLS}). 
Recall that in characteristic $p$ one
has
\[ \Gamma_{\Fp}[\alpha_2]\iso
\Fp[e_2,e_{2p},e_{2p^2},\ldots]/(e_2^p,e_{2p}^p,\ldots)
\]

It is easy to see that the map $\HH^2(\Fp,\Fp) \ra \THH^2(\Fp,\Fp)$ is
an isomorphism---the $k$-invariants in $\HH^2(\Fp,\Fp)$ classify the
two dgas $\Z/p^2$ and $\Z/p[\epsilon]/\epsilon^2$, and the
$k$-invariants in $\THH^2(\Fp,\Fp)$ classify the ring spectra
$H\Z/p^2$ and $H(\Z/p[\epsilon]/\epsilon^2)$.  So by choosing our
generators appropriately we can assume $\sigma_2$ is sent to $\alpha_2.$
We therefore have that $\sigma_2^p \in \HH^{2p}(\Fp,\Fp)$ maps to zero
in $\THH^{2p}(\Fp,\Fp)$, using the ring structure.

Let $C$ be the $\Z$-dga whose homology is $\Lambda_{\Fp}(g_{2p-2})$ and
whose nontrivial $k$-invariant is $\sigma_2^p$.  Let $D$ be the dga
$\Lambda_{\Fp}(g_{2p-2})$ with zero differential.  We know $C$ is not
quasi-isomorphic to $D$, as they have different $k$-invariants in $\HH^*$.
However, the $k$-invariants for the ring spectra $HC$ and $HD$ are the
same (and are equal to the zero element of $\THH^{2p}(\Fp,\Fp)$).
So $HC\he HD$, that is to say $C$ and $D$ are topologically
equivalent.

To be even more concrete, take $p=2$.  Then $C\he
\Z[e_1;de_1=2]/(e_1^4)$ and $D=\Lambda_{\Ft}(g_2)$. We have shown that
these are topologically equivalent, but not quasi-isomorphic.

It's worth observing that as $p$ increases the dgas produced by this
example become more complicated; to construct them explicitly one is
required to go further and further out in the resolution of $\Fp$ as a
$\Z$-dga.

\begin{remark}
In the above example one could have replaced $\sigma_2^p$ by any class
$\tau$ in the kernel of $\HH^{*}(\Fp,\Fp) \to \THH^{*}(\Fp,\Fp)$.
In particular, we could have taken $\tau$ to be any power of
$\sigma_2^p$.  It is more difficult to work out explicitly what dgas
arise from these higher $k$-invariants.
\end{remark}

\begin{remark}
Note that $\HH^*(\Fp,\Fp)$ and $\THH^*(\Fp,\Fp)$ are isomorphic as
abstract groups (equal to zero in odd dimensions, $\Fp$ in even
dimensions).  It is somewhat of a surprise that the map between them
is not an isomorphism.  The reader should take note that even without
knowledge of the ring structures on $\HH^*$ and $\THH^*$, one can
still see that $\sigma_2^p$ must map to zero.  The map $\HH^{2p} \ra
\THH^{2p}$ has the form $\Z/p \ra \Z/p$, and the $k$-invariant for the
ring spectrum $P_{2p-2}k(1)$ is certainly a non-trivial element in the
latter group (here $k(1)$ is the first Morava $K$-theory spectrum; see
\cite{A} or \cite{G} for the ring structure when $p=2$).
But this ring spectrum cannot possibly be $H(\blank)$ of any dga,
since the underlying spectrum is not Eilenberg-Mac\,Lane.  So the map
$\Z/p \ra \Z/p$ cannot be surjective---and hence not injective,
either.
\end{remark}

\subsection{Examples using $H\Z$-algebra structures}
The following examples give another approach to producing
topologically equivalent dgas.  For these examples, note that if $A$ is
a $\Z$-dga then by the \dfn{(derived) mod 2 homology ring} of $A$ we mean the
ring $H_*(A\tens^L \Z/2)$.

\begin{example}
\label{ex:tpwe1}
The problem is to give two dgas which are topologically equivalent but
not quasi-isomorphic.  This is equivalent---via the identification of
dgas and H$\Z$-algebras from~\cite{S}---to giving two $H\Z$-algebras
which are weakly equivalent as ring spectra, but not as
$H\Z$-algebras.  Consider the ring spectrum $H\Z \Smash_S H\Z/2$,
which has two obvious $H\Z$-algebra structures (from the left and
right sides of the smash).  Let $HC$ and $HD$ denote
these two different $H\Z$-algebra structures, coming from the left
and from the right respectively.   We claim that
$HC$ and $HD$ are not weakly equivalent as
$H\Z$-algebras---although obviously they are the same ring spectrum.
This will give another example of topological equivalence.

To see that $HC$ and $HD$ are distinct $H\Z$-algebras, we can give the
following argument (it is inspired by one shown to us by Bill Dwyer).
If they {\it were\/} equivalent as $H\Z$-algebras, then one would have
an equivalence of ring spectra $HC\Smash_{H\Z} E \he HD\Smash_{H\Z} E$
for any $H\Z$-algebra $E$; so there would be a resulting isomorphism of
rings $\pi_*(HC\Smash_{H\Z} E)\iso \pi_*(HD\Smash_{H\Z} E)$.  Write
$H=\HZ$ and let's choose $E=H\Z/2$.  So we will be computing the
(derived) mod $2$ homology rings of the associated dgas $C$ and $D$,
since $\pi_*(HC\Smash_{H\Z} E) \iso H_*( C \otimes^L \Z/2).$

For $HC$ we have
\[ HC\Smash_{H} E = E\Smash_H HC = E\Smash_H (H\Smash E) = E\Smash E.\]
So $\pi_*(HC\Smash_H E)=\pi_*(E\Smash E)$ and we have that the mod $2$
homology ring of $C$ is the dual Steenrod algebra.

For $HD$, however, the situation is different.  The structure map from
$H\Z$ to $HD$ factors through $\HZt$, and so $D$ is an $\F_2$-dga.  One
readily checks that for any $\F_2$-dga the mod $2$ homology has an
element of degree $1$ whose square is zero (in fact if $X$ is an
$\F_2$-dga then $H_*(X\tens^L \Z/2)\iso H_*(X) \tens_{\F_2}
\Lambda_{\F_2}(e_1)$.  We find, therefore, that $C$ and $D$ have
different mod $2$ homology rings---since the dual Steenrod algebra does
not have an element in degree $1$ squaring to zero.  So $C$ and
$D$ cannot be quasi-isomorphic.
\end{example}

\begin{example}
\label{ex:tpwe2}
In the previous example, the $\HZ$-algebras $HC$ and $HD$ are quite
big---having non-vanishing homotopy groups in infinitely many degrees.
One can obtain a smaller example by letting $HC$ and $HD$ be $\HZ
\Smash_{bo} \HZt$, with the $\HZ$-algebra structure coming from the
left and right, respectively.  One applies the same analysis as before
to see that $C$ and $D$---the associated dgas---have different mod $2$
homology rings.  One only needs to know that
$\pi_*(HC\Smash_{H} E) \iso \pi_*(E\Smash_{bo} E)
\iso A(1)_* = A_*/ (\xib_1^4, \xib_2^2, \xib_3, ...)$,
where $A_*$ is the dual Steenrod algebra and $\xib_i = c(\xi_i)$ where
$c$ is the anti-automorphism.  In particular, this ring does not have
an element of degree $1$ which squares to zero.

Note that the homotopy rings of $C$ and $D$ are just
$\Lambda_{\Ft}(f_2,g_3)$.  To see this, use our knowledge of
$\pi_*(E\Smash_{bo} E)$
together with an analysis of
the cofiber sequence of spectra
\[ H\Smash_{bo} E \llra{2} H\Smash_{bo} E \ra E \Smash_{bo} E
\]
(note that the multiplication by $2$ map is null).  This shows that
the homotopy of $C$ is $\Z/2$ in degrees $0$, $2$, $3$, and $5$---but
it does not immediately identify the class in degree $5$ as a product.
However, the mod $2$ homology of $C$ is $A(1)_*\iso
\F_2[x,y]/(x^4,y^2)$ where $|x|=1$ and $|y|=2$ (as $E\Smash_{H} HC \he
E\Smash_{bo} E$), and this ring structure is only consistent with the
degree $5$ class in $H_*(C)$ being the product of the classes in
degrees $2$ and $3$.  This argument works with the product
being taken in either order.

In fact one can determine the dgas $C$ and $D$ explicitly: $D$ is
$\Lambda_{\Ft}(g_2,h_3)$ with zero differential, and $C$ is the dga
$\Z\langle e,h;de=2, dh=0 \rangle/(e^4,h^2,eh+he)$ with $|e| = 1,
|h|=3$.  For $D$ this
follows because it is an $\Ft$-algebra, and an analysis of the
possible $k$-invariants in the Postnikov tower shows there is only one
$\Ft$-dga with the given homology ring.  For $C$, when analyzing the
Postnikov tower one finds there are two possibilities for $P_2C$ (the
two dgas given near the end of Section~\ref{ex:tpwe}).  The only one which
gives the correct mod $2$ homology of $C$ is $\Z[e; de =2]/(e^4)$.
After this stage, at each level of the Postnikov tower there is only
one possible $k$-invariant consistent with the given homology ring.

Finally, note that we cannot get a similar example by using the left
and right $H\Z$-algebra structures on $\HZ \Smash_{bu} \HZt$.  Unlike
the previous examples, these actually give weakly equivalent
$\HZ$-algebras.  To see this, apply $(\blank)\Smash_{bu} \HZt$ to the
cofiber sequence $\Sigma^2 bu \ra bu \ra \HZ$; this shows that the
homotopy ring of $\HZ\Smash_{bu}\HZt$ is $\Lambda_{\Ft}(e_3)$.  We
showed in Example~\ref{ex:exterior} that there is only one homotopy
type for dgas with this homology ring.
\end{example}

\subsection{Topological equivalence over fields}
\label{se:fields}

For the above examples of nontrivial topological equivalence, it has
been important in each case that we were dealing with dgas over $\Z$
whose zeroth homology is $\Z/p$.  In each example one of the dgas
involved a class in degree $1$ whose differential is $p$ times the
unit.  This leads to the following two questions:

\smallskip

\noindent
{\it Question 1\/}: Do there exists two dgas over $\Z$ which have
torsion-free homology, and which are topologically equivalent but not
quasi-isomorphic?

\smallskip\noindent
{\it Question 2\/}: Let $F$ be a field.  Do there exist two
$F$-dgas which are topologically equivalent but not quasi-isomorphic
(as $F$-dgas)?

\smallskip

So far we have been unable to answer these questions except in one
simple case.  We can show that if two $\Q$-dgas are topologically
equivalent then they must actually be quasi-isomorphic:

\begin{proof}[Proof of Proposition~\ref{prop-Q-same}]
The key here is that for any ring spectrum $R$ with rational homotopy,
the map $\eta\colon R=R\Smash_S S \to R \sm_S H\bQ$ is a weak
equivalence.  One way to see this is to note that $\bQ$ is flat over
$\pi_*^s$, as it is the localization of $\pi_* S$ at the set of
nonzero integers.  So the spectral sequence calculating the homotopy
of $X \sm_S H\bQ$ from~\cite[IV.4.2]{ekmm} collapses, showing that
$\eta$ induces isomorphisms on homotopy.

If $A$ is a
\Qdga then the above shows that $\eta_A\mc HA \to HA \sm_S H\bQ$ is a weak
equivalence of ring spectra.

Assume $C$ and $D$ are two topologically equivalent \QdgasP.  Then
$HC$ and $HD$ are weakly equivalent as $S$-algebras.  It follows that
$HC \sm_S H\bQ$ and $HD \sm_S H\bQ$ are equivalent as $H\bQ$-algebras.
If $\eta_C$ and $\eta_D$ were $H\bQ$-algebra maps we could conclude
that $HC$ and $HD$ were equivalent as $H\bQ$-algebras, and hence that
$C$ and $D$ are quasi-isomorphic \QdgasP.  Unfortunately, the claim
that $\eta_C$ and $\eta_D$ are $H\bQ$-algebra maps is not so clear.
Instead, we will use the
map $\psi_C \mc HC \sm_S H\bQ \to HC \sm_{H\bQ} H\bQ \iso HC$.  This
is a map of $H\bQ$-algebras.  Note that $\psi_C \eta_C$ is the
identity, and so $\psi_C$ is also a weak equivalence.

Using  $\psi_C$ and $\psi_D$ we obtain a zig-zag of weak equivalences of
$H\bQ$-algebras $HC \bwe HC\Smash_S H\bQ \he HD\Smash_S H\bQ \we HD$.
So $C$ and $D$ are quasi-isomorphic as $\Q$-dgas.
\end{proof}


\section{Homotopy endomorphism spectra and dgas}
\label{se:hoend}

The next main goal is the proof of our Tilting Theorem (\ref{thm-tilting}).
The portion of the proof requiring the most technical difficulty
involves keeping track of information preserved by a zig-zag of
Quillen equivalences.  The machinery needed to handle this is
developed in \cite{hend} and \cite{additive}.
The present section summarizes what we need.

\medskip

A model category is called \dfn{combinatorial} if it is
cofibrantly-generated and the underlying category is locally
presentable.  See \cite{comb} for more information.  The categories of
modules over a dga and modules over a symmetric ring spectrum are both
combinatorial model categories.

If $\cM$ is a combinatorial, stable model category then \cite{hend}
explains how to associate to any object $X\in \cM$ a \dfn{homotopy
endomorphism ring spectrum} $\hEnd(X)$.  This should really be
regarded as an isomorphism class in $\Ho(\RingSpectra)$, but we will
usually act as if a specific representative has been chosen.

\begin{prop}
\cite[Thm. 1.4]{hend} 
\label{pr:hoend}
Let $\cM$ and $\cN$ be combinatorial, stable model categories.
Suppose that $\cM$ and $\cN$ are connected by a zig-zag of Quillen
equivalences (in which no assumptions are placed on the intermediate
model categories in the zig-zag).  Suppose that $X\in \cM$ and $Y\in
\cN$ correspond under the derived equivalence of homotopy categories.
Then $\hEnd(X)$ and $\hEnd(Y)$ are weakly equivalent ring spectra.
\end{prop}

Recall that a category is said to be {\it additive\/} if the Hom-sets
have natural structures of abelian groups, the composition is
bilinear, and the category has finite coproducts.  See \cite[Section
VIII.2]{MacL}.  Such a category is necessarily pointed: the empty
coproduct is an initial object, which is also a terminal object by
\cite[VIII.2.1]{MacL}.

By an \dfn{additive model category} we mean a model category whose
underlying category is additive and where the additive structure
interacts well with the `higher homotopies'.  See \cite[Section
2]{additive} for a precise definition, which involves the use of
cosimplicial resolutions.  If $R$ is a dga, the model category $\Mod
R$ of differential graded $R$-modules is one example of an additive
model category; this example is discussed in more detail at the
beginning of the next section.

Note that if $L\colon \cM \adjoint \cN\colon R$ is a Quillen pair
where $\cM$ and $\cN$ are additive, then both $L$ and $R$ are additive
functors---this is  because they preserve direct sums (equivalently,
direct products) since $L$ preserves colimits and $R$ preserves
limits.

If $\cM$ and $\cN$ are two additive model categories, we say they are
\dfn{additively Quillen equivalent} if they can be connected by a
zig-zag of Quillen equivalences where all the intermediate categories
are additive.
As remarked in the introduction, it is possible for two additive model
categories to be Quillen equivalent but not additively Quillen
equivalent.  We'll give an example in Section~\ref{se:QEnotAQE} below.

\subsection{$\ChZ$ enrichments}
Recall that $\ChZ$ denotes
the model category of chain complexes of abelian groups, with its
projective model structure.  A $\ChZ$-model category is a model
category with compatible tensors, cotensors, and enrichments over
$\ChZ$; see \cite[4.2.18]{hovey} or \cite[Appendix A]{hend}.  For $X,
Y$ in $\cM$, we denote the function object in $\ChZ$ by $\Fhomu(X,Y)$.

Note that a pointed $\ChZ$-model category is automatically an additive
and stable model category.  The additivity of the underlying
category, for instance, follows from the
isomorphisms $\cM(X,Y)\iso \cM(X\tens \Z,Y)\iso
\Ch(\Z,\Fhomu(X,Y))$ and the fact that the last object has a
natural structure of abelian group.   This additive structure is also
compatible with the higher homotopies; see~\cite[2.9]{additive}. 
The stability follows directly
from the stability of $\ChZ$; cf.~\cite[3.5.2]{SS2} 
or~\cite[3.2]{GS}. 

We will need the following result, which is proven in \cite{additive},
connecting $\Ch$-enrichments to homotopy endomorphism spectra:

\begin{prop}
\label{pr:hoend2}
Let $\cM$ be a combinatorial $\ChZ$-model category such that $\cM$
has a generating set of compact objects
(cf. Definition~\ref{def:CompGen}).
Let $X\in \cM$ be a cofibrant-fibrant object.  Then
\begin{enumerate}[(a)]
\item $\hEnd(X)$ is weakly equivalent to the Eilenberg-Mac\,Lane ring
spectrum associated to the dga $\Fhomu(X,X)$.
\item Suppose $\cN$ is another combinatorial, $\ChZ$-model category
with a generating set of compact objects, and that $\cM$ and $\cN$ are
connected by a zig-zag of additive Quillen equivalences.  Let $Y$ be a
cofibrant-fibrant object corresponding to $X$ under the derived
equivalence of homotopy categories.  Then the dgas $\Fhomu(X,X)$ and
$\NFhomu(Y,Y)$ are quasi-isomorphic.
\end{enumerate}
\end{prop}

Part (a) of the above result follows directly from \cite[1.4,
1.6]{additive}. 
Likewise, part (b) follows directly from \cite[1.3, 1.6]{additive}. 


\section{Tilting theory}
\label{se:tilting}

This section addresses the following question: given two dgas $C$ and
$D$, when are the model categories of (differential graded) $C$-modules and
$D$-modules
Quillen equivalent?  We give a complete answer in terms of topological
tilting theory, making use of topological equivalence.  There is also
the associated question of when the two module categories are {\it
additively\/} Quillen equivalent, which can be answered with a
completely algebraic version of tilting theory.

\medskip

If $C$ is a dga, let $\rMod C$ be the category of (right) differential
graded $C$-modules.  This has a model structure lifted from $\Ch$,
in which a map is a fibration or weak equivalence if and only if the
underlying map of chain complexes is so; see~\cite[4.1(1)]{SS1}.
We let $\cD(C)$ denote the homotopy category of $\rMod C$, and call
this the \dfn{derived category} of $C$.

The category of $C$-modules is enriched over $\Ch$: for $X, Y$ in
$\rMod C$, let $\Hom_C(X,Y)$ be the chain complex which in degree $n$
consists of $C$-module homomorphisms of degree $n$ on the underlying graded
objects (ignoring the differential).  The differential on $\Hom_C(X,Y)$
is then defined so that the chain maps are the cycles.
See~\cite[4.2.13]{hovey}.

\begin{defn}\label{def:CompGen}
Let $\cT$ be a triangulated category with infinite
coproducts.
\begin{enumerate}[(a)]
\item
An object $P\in \cT$ is called \dfn{compact} if
$\oplus_\alpha \cT(P,X_\alpha) \ra \cT(P,\oplus_\alpha
X_\alpha)$ is an isomorphism for every set of objects $\{X_\alpha\}$;
\item A set of objects $S\subseteq \cT$ is a \dfn{generating set} if
the only full triangulated subcategory of $\cT$ which contains $S$
and is closed under arbitrary coproducts is $\cT$ itself.
If $S$ is a singleton set $\{P\}$ we say that $P$ is a
\dfn{generator}.
\end{enumerate}
\end{defn}

An object $X$ in a stable model category $\cM$ is called
compact if it is a compact object of the associated homotopy
category.  Likewise for the notion of a generating set.

We can now state the main result of this section:

\begin{thm}[Tilting Theorem]\label{thm-tilting}
Let $C$ and $D$ be two dgas.
\begin{enumerate}[(a)]
\item
The model categories of $C$-modules and $D$-modules are Quillen
equivalent if and only if there is a cofibrant and fibrant
representative $P$ of a compact generator in $\cD(C)$
such that $\Hom_C(P,P)$ is topologically equivalent to $D$.
\item
The model categories of $C$-modules and $D$-modules are {\em additively}
Quillen equivalent if and only if there is
a cofibrant and fibrant representative $P$ of
a compact generator in
$\cD(C)$ such that $\Hom_C(P,P)$ is quasi-isomorphic to $D$.
\end{enumerate}
\end{thm}

Before proving this theorem, we need to recall a few results on ring
spectra and their module categories.  If $R$ is a ring spectrum, we
again let $\rMod R$ denote the category of right $R$-modules equipped
with the model structure of \cite[4.1(1)]{SS1}. 
  This model
category is enriched over symmetric spectra: for $X,Y\in \rMod R$ we
let $\Hom_R(X,Y)$ denote the symmetric spectrum mapping object.

\begin{prop}
\label{prop-6.3}
\begin{enumerate}[(a)]
\item
\cite[4.1.2]{SS2}\label{prop-one-gen}
Let $R$ be a ring spectrum, and let $P$ be a cofibrant and fibrant
representative of a compact generator of the
homotopy category of $R$-modules.  Then there is a Quillen equivalence
$\rMod \Hom_R(P,P) \ra \rMod R$.
\item\cite[3.4]{GS}\label{prop-add-gen} 
Let $C$ be a dga, and let $P$ be a cofibrant and fibrant representative
of a compact generator of the
homotopy category of $C$-modules.  Then there is a Quillen equivalence
$\rMod \Hom_C(P,P) \ra \rMod C$.
\end{enumerate}
\end{prop}

In~\cite[3.4]{GS} this is actually proved over $\Q$, but the same
proofs work over $\Z$.

\begin{prop}
\label{prop-6.4}
\begin{enumerate}[(a)]
\item
\cite[5.4.5]{hss}\label{prop-ho-inv}
If $R\ra T$ is a weak equivalence of ring spectra, then there is an
induced Quillen equivalence $\rMod R\ra \rMod T$.
\item\cite[4.3]{SS1}\label{prop-dga-inv}
If $C \ra D$ is a quasi-isomorphism of dgas, then
there is an induced  Quillen equivalence $\rMod C\ra \rMod D$.
\end{enumerate}
\end{prop}

Now we can prove the Tilting Theorem:

\begin{proof}[Proof of Theorem~\ref{thm-tilting}]
For part (a), note first that if $\rMod C$ and $\rMod D$ are Quillen
equivalent then $\cD(C)$ and $\cD(D)$ are triangulated equivalent.
Let $P$ be the image of $D$ in $\cD(C)$.  Since $D$ is a compact
generator of $\cD(D)$, its image $P$ is a compact generator of
$\cD(C)$.  But $\rMod C$ and $\rMod D$ are combinatorial, stable model
categories, so by Proposition~\ref{pr:hoend} we know $\hEnd(P)$ and
$\hEnd(D)$ are weakly equivalent ring spectra.  The two module
categories are also $\Ch$-categories, so if follows from
Proposition~\ref{pr:hoend2}(a) that
\[ \hEnd(P)\he H(\Hom_C(P,P)) \qquad\text{and}\qquad
\hEnd(D)\he H(\Hom_D(D,D)).
\]
Since $\Hom_D(D,D)$ is isomorphic to $D$, we have that
$\Hom_C(P,P)$ and $D$ are topologically equivalent.

For the other direction, we are given that $H\Hom_C(P,P)$ and $HD$ are
weakly equivalent as ring spectra.  Thus $H\Hom_C(P,P)$-modules and
$HD$-modules are Quillen equivalent by
Proposition~\ref{prop-6.4}(\ref{prop-ho-inv}).  By~\cite[2.8]{S}, the
model category of $\Hom_C(P,P)$-modules is Quillen equivalent to
$H\Hom_C(P,P)$-modules and similarly for $D$-modules and $HD$-modules.  So
$\rMod D$ and $\rMod \Hom_C(P,P)$ are Quillen equivalent.
Proposition~\ref{prop-6.3}(\ref{prop-one-gen}) then finishes the
string of Quillen equivalences by showing that $\rMod C$ and
$\rMod \Hom_C(P,P)$ are Quillen equivalent.

Now we turn to part (b) of the theorem.  If the categories of
$C$-modules and $D$-modules are additively Quillen equivalent then the
image of $D$ in $\cD(C)$ is a compact generator, just as in part (a).
By Proposition~\ref{pr:hoend2}(b), $\Hom_D(D,D)$ is quasi-isomorphic to
$\Hom_C(P,P)$.  We have already remarked that $D$ and $\Hom_D(D,D)$ are
isomorphic, so $D$ is quasi-isomorphic to $\Hom_C(P,P)$.

For the other direction, suppose given a compact generator $P$ in $\cD(C)$.
Proposition~\ref{prop-6.3}(\ref{prop-add-gen}) shows that $\rMod C$ is
additively Quillen equivalent to $\rMod \Hom_C(P,P)$, and
Proposition~\ref{prop-6.4}(b) shows that
$\rMod \Hom_C(P,P)$ is additively Quillen equivalent to $\rMod D$.
\end{proof}

\begin{remark}
\label{re:equivalences}
When $R$ and $S$ are rings, the following statements are equivalent:
\begin{enumerate}[(1)]
\item $\cD(R)$ is triangulated-equivalent to $\cD(S)$;
\item There is a cofibrant and fibrant representative $P$ of a compact
generator in $\cD(S)$ such that
the dga $\Hom_S(P,P)$ is quasi-isomorphic to $R$;
\item $\Ch_R$ and $\Ch_S$ are Quillen equivalent model categories.
\end{enumerate}
The equivalence of (1) and (2) was established by Rickard \cite{R},
and the equivalence with (3) was explicitly noted in
\cite[2.6]{dereq} 
 (this reference only discusses {\it pointed\/}
Quillen equivalences, but that restriction can be removed using
\cite[5.5(b)]{hend}: 
 if $\Ch_R$ and $\Ch_S$ are Quillen
equivalent, then they are Quillen equivalent through a zig-zag of
pointed model categories).

When $R$ and $S$ are dgas---rather than just rings---the situation
changes somewhat.  Theorem~\ref{thm-tilting} establishes that the
analogs of (2) and (3) are still equivalent, where in (2)
``quasi-isomorphic'' is replaced by ``topologically equivalent''.  And
(3) certainly implies (1).  But the implication $(1)\Rightarrow (2)$
is not true.  One counterexample is $R=\Z\langle e, x, x^{-1} ;
de=p, dx=0\rangle/ (ex+xe=x^2)$ with $|e| = 1, |x| = 1$ and 
$S=H_*(R)=\Z/p[x_1,x_1^{-1}; dx_1=0]$.  The
verification that this is indeed a counterexample will be taken up in
the paper~\cite{stabmod}.  The dgas $R$ and $S$ arise in connection
with the
stable module categories discussed in~\cite{marco}.
\end{remark}

\subsection{Relation to the theory of dg-categories}
\label{se:dgcats}
For dgas over the ground ring $\Z$ there are several different
notions of equivalence which one might consider.  Here are
five of them:
\begin{enumerate}[(1)]
\item Quasi-isomorphism;
\item Topological equivalence;
\item Additive Morita equivalence, meaning that the model categories
  of $R$-modules and $T$-modules are additively Quillen equivalent;
\item Topological Morita equivalence, meaning that the model
  categories of $R$-modules and
  $T$-modules are Quillen equivalent;
\item dg-equivalence, meaning that the dg-categories of 
  cofibrant-fibrant objects in $R$-modules
  and $T$-modules are {\it quasi-equivalent\/} (explored in papers
  such as \cite{keller.dg,keller,tabuada,toen}).    
\end{enumerate}
Here ``dg-category'' just means a category enriched over the symmetric
monoidal category of chain complexes over $k$.  A dg-category $\cG$
has a corresponding `homotopy category' $\ho(\cG)$ with the same
objects as $\cG$, and where the morphisms from $x$ to $y$ are the
zeroth homology of the chain complex $\cG(x,y)$.  A quasi-equivalence
of dg-categories is a functor $f\colon \cC \to \cD$ which induces
quasi-isomorphisms on all of the morphism complexes $\cC(x,y) \to
\cD(x,y)$, and which induces an equivalence of homotopy categories.
As usual, two dg-categories are said to be quasi-equivalent if they can be
connected by a zig-zag of quasi-equivalences.

One has the implications
\[ \xymatrix{
(1) \ar@{=>}[r]\ar@{=>}[d] & (3) \ar@{=>}[d]\ar@{=>}[r] & (5) \\
(2) \ar@{=>}[r] & (4),
}
\] 
all of these being elementary except $(3)\Rightarrow(5)$. 
Theorem~\ref{thm-tilting} (and the results of \cite{additive}) show that
any additive Morita equivalence is given by tilting functors
as in Propositions~\ref{prop-6.3}(b) and~\ref{prop-6.4}(b).
These functors preserve the dg-enrichments and
thus induce a quasi-equivalence on the
subcategories of cofibrant-fibrant objects. 

Except possibly for
$(3)\Rightarrow (5)$, none of the implications is reversible.  For the
vertical implications this is justified in the present paper.  For the
horizontal implications one can give counterexamples just using
ordinary rings: this is classical tilting theory as in
\cite{R}, together with the equivalences discussed in 
Remark~\ref{re:equivalences}.
We have not been able to decide whether $(3)\Rightarrow(5)$ is
reversible.  This is perhaps an interesting question.

In \cite[Section 3.8]{keller}, 
Keller defines two dgas $R$ and $T$ to be \dfn{dg Morita equivalent}
if there is an equivalence $\ho(\Mod R)\he \ho(\Mod T)$ given by a
composition of tensor functors and their inverses.  It turns out that
this is the same as the notion of additive Morita equivalence from
(3).  The proof comes from comparing \cite[8.2]{keller.dg} 
(quoted in~\cite[3.11]{keller}), which characterizes $dg$-Morita
equivalence in terms of tilting theory, to our
Theorem~\ref{thm-tilting} above.  We include this result for completeness.

\begin{prop}
Two $\Z$-dgas $R$ and $T$ are dg Morita equivalent if and only if the
model categories $\Mod R$ and $\Mod T$ are additively Quillen equivalent.
\end{prop}

\begin{proof}
If $\Mod R$ and $\Mod T$ are additively Quillen equivalent, then
Theorem~\ref{thm-tilting} says that there is a cofibrant, compact
generator $G$ of $\Mod T$ whose endomorphism dga is quasi-isomorphic to
$R$.  Thus, condition (2) of \cite[Thm. 3.11]{keller} is satisfied,
and so $R$ and $T$ are dg Morita equivalent.

Suppose conversely that $\Mod R$ and $\Mod T$ are dg Morita equivalent.
Let $\bar{R}$ denote the dg-category with one object whose
endomorphism dga is $R$.
By \cite[Thm. 3.11]{keller}, there is a subcategory $\cG$ of
$\Mod T$ whose objects are a generating set of compact, cofibrant
objects, such that the dg-category determined by $\cG$ is
quasi-equivalent to $\bar{R}$.  Here the dg-structure on $\cG$ is
inherited from that of $\Mod T$.  

Recall that a dg-category  has a corresponding `homotopy
category', and quasi-equivalent dg-categories have equivalent homotopy
categories.  In our case $\ho(\cG)$ is just the subcategory of
$\ho(\Mod T)$ determined by the objects of $\cG$.  

Since $\cG$ is quasi-equivalent to $\bar{R}$, the
homotopy category of $\cG$ is equivalent to the homotopy category of
$\bar{R}$.  Since the latter has exactly one object, it follows that
all the objects of $\ho(\cG)$ are isomorphic.  That is to say, all the
objects of $\cG$ are isomorphic in $\ho(\Mod T)$.  Since the objects of
$\cG$ were a generating set, it follows that every object of $\cG$ is
itself a generator for $\Mod T$.

A model category structure on $dg$-categories is constructed in
\cite{tabuada}, and later used in \cite{toen}.  The weak equivalences
are the quasi-equivalences.  As remarked in \cite[2.3]{toen}, there is a 
cofibrant-replacement $Q\bar{R}\we \bar{R}$ where the map on objects
is the identity.  Note that the endomorphism dga of the unique object
of $Q\bar{R}$ is quasi-isomorphic to $R$.

Since $\bar{R}$ is quasi-equivalent to $\cG$ and all $dg$-categories
are fibrant, 
there is a quasi-equivalence $Q\bar{R} \ra \cG$.  It follows that for
the object of $\cG$ in the image of this functor, the endomorphism dga
is quasi-isomorphic to $R$.  Hence we have a cofibrant, compact
generator of $\Mod T$ whose endomorphism dga is quasi-isomorphic to
$R$.  By Theorem~\ref{thm-tilting} above, $\Mod R$ and $\Mod T$ are
additively Quillen equivalent.
\end{proof}

\begin{remark}
Let $k$ be a commutative ring.  Keller \cite{keller} actually
defines a notion of $k$-linear $dg$ Morita equivalence for $k$-dgas.
One could possibly revise the notion of an additive model
category, instead define a `k-linear model category' together with
the notion of `k-linear Quillen equivalence', and finally prove a
$k$-linear analog of our additive tilting theorem.  We have not
pursued this, however.
\end{remark}

\section{A model category example}
\label{se:QEnotAQE}
In this section we give an example of two additive model categories
which are Quillen equivalent but not additively Quillen equivalent.

\medskip

Let $C$ and $D$ be the dgas $\Z[e_1;de_1=2]/(e_1^4)$ and
$\Lambda_{\Z/2}(g_2)$.  We have already seen in
Section~\ref{ex:tpwe} that these dgas are topologically equivalent.
Therefore $\rMod C$ and $\rMod D$ are Quillen equivalent model
categories (by Theorem~\ref{thm-tilting}(a), for instance).
We claim that they are not {\it additively\/} Quillen
equivalent, however.  If they were, then by Theorem~\ref{thm-tilting}
there would be a compact generator $P$ in $\rMod D$ such that the dga
$\Hom_D(P,P)$ is quasi-isomorphic to $C$.  But since everything in
$\rMod D$ is a $\Z/2$-module, $\Hom_D(P,P)$ is a $\Z/2$-module as
well.  We will be done if we can show that $C$ is not quasi-isomorphic
to any dga defined over $\Z/2$.

Assume $V$ is a $\Z/2$-dga, and $C$ is quasi-isomorphic to $V$.  Let
$Q \ra C$ be the cofibrant-replacement for $C$ constructed as in
Example~\ref{ex:cofibrep}.  Our assumption implies that there is a
weak equivalence $Q\ra V$.  The map $Q_0 \ra V_0$ is completely
determined, since $Q_0=\Z$ and the unit must map to the unit.  The map
$Q_1 \ra V_1$ must send $e$ to an element $E\in V_1$ such that $dE=0$
(using that $2V_0=0$ and $de=2$).  But $H_1(V)=H_1(C)=0$, and so
$E=dX$ for some $X\in V_2$.  Note that we then have $d(EX)=-E^2$ by
the Leibniz rule.  However, the generator of $H_2(Q)$ is $e^2$, and we
have just seen that the image of $e^2$ is zero in homology (since
$E^2$ is a boundary).  This contradicts the map $H_2(Q)\ra H_2(V)$
being an isomorphism.  This completes our example.


\bibliographystyle{amsalpha}

\end{document}